\documentclass[12pt]{article}

\usepackage[text={400pt,660pt},centering]{geometry}

\usepackage[english]{babel}
\usepackage[latin1]{inputenc}
\usepackage[T1]{fontenc}
\usepackage{geometry}
\usepackage{amsmath}
\usepackage{amssymb}
\usepackage{amsthm}
\usepackage{mathtools}
\usepackage{bbm}
\usepackage{seqsplit}
\usepackage{xstring}
\usepackage{graphicx}
\usepackage{lmodern}
\usepackage{mathrsfs}

\usepackage{csquotes}
\usepackage[all]{xy}
\usepackage{picture}
\usepackage{wrapfig}
\usepackage{float}
\usepackage{bm}
\usepackage[hyphens]{url}
\usepackage{enumerate}
\usepackage{hyperref}
\hypersetup{pdfborder={0 0 0},
	colorlinks=true,
	linkcolor=blue,
	citecolor=red,
	urlcolor=black}
\usepackage[marginpar]{todo}
\usepackage{dsfont}
\newcommand{\pr}{\mathbb{P}}
\newcommand{\E}{\mathbb{E}}

\newcommand{\lto}{\longrightarrow}
\newcommand{\impl}{\Longrightarrow}

\newcommand{\Var}{\mbox{Var}}
\newcommand{\Cov}{\mbox{Cov}}
\allowdisplaybreaks

\newcommand{\pf}[1]{\begin{proof}{#1}\end{proof}}
\newenvironment{myproof}{\paragraph{\it{Proof of Lemma \ref{branchbound}}}}{\hfill$\square$}

\newenvironment{myproof61}{\paragraph{\it{Proof of Lemma \ref{lem61}}}}{\hfill$\square$}

\newenvironment{myproof62}{\paragraph{\it{Proof of Lemma \ref{lem62}}}}{\hfill$\square$}

\newtheorem{lem}{Lemma}[section]
\newcommand{\lm}[1]{\begin{lem}{#1}\end{lem}}
\newtheorem{theo}{Theorem}[section]
\newcommand{\thm}[1]{\begin{theo}{#1}\end{theo}}
\newtheorem{prp}{Proposition}[section]
\newcommand{\prpo}[1]{\begin{prp}{#1}\end{prp}}
\newtheorem{corr}{Corollary}[section]
\newcommand{\corl}[1]{\begin{corr}{#1}\end{corr}}
\newtheorem{rem}[theo]{Remark}
\newcommand{\rema}[1]{\begin{rem}{\ #1}\end{rem}}
\newtheorem{rems}[theo]{Remarks}

\author{}
\date{\today}

\title{Flooding and Diameter in General Weighted Random Graphs}
\author{Thomas Mountford\thanks{Institut de Math\'ematiques, Ecole Polytechnique F\'ed\'erale de Lausanne, 1015 Lausanne, Switzerland. Email: thomas.mountford@epfl.ch, jacques.saliba@epfl.ch}\and Jacques Saliba \footnotemark[1]}

\begin{document}

%

\maketitle



\begin{abstract}
\noindent We study in this paper, the first passage percolation on a random graph model, the configuration model. We first introduce, the notions of weighted diameter, which is the maximum of the weighted lengths of all optimal paths between any two vertices in the graph, and the flooding time, which represents the time  (weighted length) needed to reach all the vertices in the graph starting from a uniformly chosen vertex.
Our result consists of describing the asymptotic behavior of the diameter and the flooding time, as the number of vertices $n$ tends to infinity, in the case where the weight distribution $G$ has an exponential tail behavior, and proving that this category of distributions is the largest possible for which the asymptotic behavior holds.
\end{abstract}




\section{\textbf{Introduction}}
\noindent Many random graph models have been developed in the last decades in order to describe real world complex systems such as as social networks and the Internet. Given a connected graph with $n$ nodes, we assign positive random weights to the edges that represent the cost, the transmission information time, or the infection time for example (in an epidemic model) among the vertices. We typically assume that these weights are i.i.d.. The optimal path between two uniformly chosen vertices $u$ and $v$ is the path between them with the minimal edge weights sum. More precisely, writing $X_e\sim G$ for an edge $e$ and for a continuous distribution $G$, and writing $\Gamma_{uv}$ the set of all paths between $u$ and $v$, the weighted length $L_n=L_n(u,v)$ of the optimal path between $u$ and $v$ is given by
\begin{equation*}
L_n=\underset{\pi\in \Gamma_{u,v}}{\min} \sum_{e\in \pi} X_e.
\end{equation*}
\noindent So $L_n$ can be viewed as the infection time of the vertex $v$ knowing that $u$ is infected (or vice versa) in a network epidemic model. The \textit{diameter} of the resulting graph will be the maximum of these optimal paths for any randomly chosen pair of vertices, and the \textit{flooding} with respect to a vertex $u$ is the maximal time that we need to spend to reach all the vertices in the graph starting from $u$. Again, we use  \textit{First Passage Percolation} techniques in order to describe the asymptotic behavior of the diameter and the flooding in the weighted \textit{Configuration Model}, a random graph with prescribed degrees; F.P.P can describe how a fluid spreads in a medium. Several authors studied the asymptotic behavior of the diameter for a non-weighted random graph, as Fernholz and Ramachandran in \cite{fernholz2007diameter} and van der Hofstad, Hooghiemstra and Znamenski in \cite{hofstad2007phase}.\\
Bhamidi, van der Hofstad and Hooghiemstra obtained the asymptotic distributions of the typical weight between two randomly chosen vertices and of the hopcount, which is the number of edges in the optimal path, in the exponential weight case at first \cite{bhamidi2010first} and in the general case \cite{bhamidi2017universality}. Amini, Lelarge and Draief obtained a law of large numbers of the diameter and the flooding in the configuration Model with exponential edge weights \cite{flood},\cite{amini2015diameter}. We give, in this paper, a generalization of their results to all edge weight distributions having a certain exponential tail behavior.  

\section{\textbf{Definitions and notations}}
\noindent We first recall the well known \textit{Configuration Model} described in details in \cite{van2016random} and \cite{bhamidi2017universality}. Given an integer $n$ and a sequence $\mathbf{d}:=(d_i^n)_{i=1}^n$ of non-negative integers such that $\sum_{i=1}^n d_i^n$ is even, the \textit{Configuration model} on $n$ vertices is constructed as follows:\\
We start with $n$ vertices numbered from 1 to $n$, and we assign $d_i^n$ half-edges to the $i$th vertex. The random graph $CM_n(\mathbf{d})$ is obtained by randomly choosing pairs of half-edges to form edges between the two corresponding vertices. 
 Let $F_n$ be the cumulative distribution of the degree of a randomly chosen vertex, denoted by $D_n$, that is
\begin{equation*}
F_n(x)=\frac{1}{n}\sum_{i=1}^n \mathbbm 1_{\{d_i^n\leq x\}}.
\end{equation*}
We let $V^n$ denote the set of vertices $\{1,2, \cdots n\}$ and $l_n:=\sum_{i\in V^n} d_i$ be the total degree of the graph.
We assume that there exists a distribution $\textbf{p}=(p_k)_{k\geq 0}$ such that $\mathbf{d}$ and \textbf{p} satisfy the following regularity conditions, as in \cite{flood}:\\
\ \\
\phantomsection
\textit{\underline{Condition 1:}}
\begin{enumerate}[a)] \label{enum}
\item $\frac{\#\{i\ | \ d_i^n=r\}}{n}\to p_r\  \forall r\geq0, n\to \infty$,
\item $\min_{i=1,\cdots, n}d_i^n:=d_{\min}\geq 3$ and $p_{d_{\min}}>0$,
\item $\limsup_{n \rightarrow \infty } \frac{1}{n} \sum_i (d_i^n) ^ {2+\delta}<\infty$ for a certain $\delta>0$.
\end{enumerate}
\rema{\label{con-d}
 Condition $(c)$ above ensures the convergence of first and second moments of $D_n$ to the respective moments of $D$ (a random variable distributed according to $(p_r)_{r\geq d_{\min}}$). Moreover, it gives an upper bound for the maximal degree $\Delta_n$ of the graph constructed on $n$ vertices. Indeed,
this condition is equivalent to $\sup_n \E[D_n^{2+\delta}\log D_n]<\infty$ and so $\E[D_n^2]$ is uniformly upper bounded. By the uniform integrability of $D_n^2$, we get
\begin{equation*}
\sum_{r\geq 3} r^2p_r^{(n)}\to\sum_{r\geq 3}r^2p_r,
\end{equation*}
where $p_r^{(n)}= \frac{\#\{i\ | \ d_i^n=r\}}{n}$. The argument is similar for $\E[D_n]\to \E[D]$.\\
On the other hand, by writing $\Delta_n$ for the maximal degree in $CM_n(\textbf{d})$, we have
\begin{equation}\label{knr}
\Delta_n^{2+\delta}=o(n)\Longrightarrow \Delta_n=o\left(n^{1/{(2+\delta)}}\right)\Longrightarrow \Delta_n=o\left(\sqrt{n/\log n}\right).
\end{equation}
 
}

Under these conditions the resulting random graph may have loops or  multi-edges but we will see that locally, the random graph will not have either and will look like a random tree. This will be detailed as a coupling argument in sections \ref{upperbound} and \ref{slb}, based on \cite{flood} and \cite{bhamidi2017universality}. In fact for a vertex $v$ picked at random among $\{1,2 \cdots n\}$, the number of vertices at (graphical) distance $r$ from $v$ will tend in distribution, as $n$ tends to infinity, to that of an inhomogeneous branching process which for generation 1 has distribution $\textbf{p}=(p_k)_{k\geq 0}$ for the number of offspring but thereafter has 
the \textquote{size biased} distribution
\begin{equation}\label{bias}
\widehat{\textbf{p}}=(\widehat{p_k})_{k\geq 0}  \mbox{ for } \widehat{p_k} = \frac{(k+1) p_{k+1}}{m}
\end{equation}
for the number of offspring.
The assumption (c) in \hyperref[enum]{Condition 1} guarantees that the distribution $\widehat{\textbf{p}}=(\widehat{p_k})_{k\geq 0}$ has finite mean (which we denote by $\nu$).  Note that $\nu $ is greater than $d_{\min} -1 \geq 2$.

We recall the \textit{Malthusian parameter} $\alpha$ corresponding to the rate at which a continuous time branching process grows, with splitting law $\widehat{\textbf{p}}=(\widehat{p_k})_{k\geq 0}$ and lifetimes distributed as $G$. It is the unique positive real number satisfying 
\begin{equation}\label{malt}
\nu  \int_0 ^ \infty e^{- \alpha t} dG(t) \ = \ 1.
\end{equation}
\noindent The population of the branching process will grow at rate $\alpha $.\\
The following distribution, that tends to the size-biased distribution as $n\to \infty$, will be used for the upper bound of the diameter:
\begin{equation*}
p_k^n:=\frac{k+1}{l_n}\sum_{i=1}^n \mathbbm 1_{\{d_i=k+1\}}.
\end{equation*}
We denote $\nu_n$ its mean and $\alpha_n$ its corresponding Malthusian parameter. It is easy to see that $\nu_n\to \nu$, and so we have that $\alpha_n\to \alpha$ as $n\to \infty$ using the fact that $l_n/n=\frac{1}{n}\sum_{i=1}^n d_i^{(n)}$ tends to $m$ by \hyperref[enum]{Condition 1}.
\\ We give $i.i.d$ positive random weights for the edges following a continuous law $G$ that has an exponential tail behavior, that is:
\begin{equation}\label{*}
\underset{x\to\infty}{\lim} \frac{-\log\overline{G}(x)}{x}=c\in (0,\infty),
\end{equation} 
where $\overline{G}(x):=1-G(x)$.\\
We write $dist_w(a,b)$ for the sum of the weights along the optimal path between $a$ and $b$, the weights being $i.i.d$ according to a continuous law $G$ satisfying \eqref{*}.  We define the weighted diameter and the weighted flooding time of $CM_n(\mathbf{d})$ as
\begin{align*}
diam(CM_n(\mathbf{d}))&:=\max\{dist_w(a,b),a,b\in V\},\\
flood(CM_n(\mathbf{d}))&:=\max\{dist_w(a,b),b\in V\},
\end{align*}
where $V$ is the set of vertices of $CM_n(\mathbf{d})$, and where the vertex $a$ in the $flood$ is chosen uniformly at random in the flooding definition.\\
For the sequel of this paper, we say that an event $A_n$ holds with high probability (w.h.p.) when $\pr(A_n)\to 1$ as $n\to \infty$.\\
 The same methods used in this paper joint with complementary arguments can be used to derive the general case (where $d_{\min}\geq 1)$, similarly to \cite{amini2015diameter}. 
\subsection{Exploration process}\label{expro}
We use, instead of constructing the random graph and then looking for the optimal path between two vertices, a coupling argument as in \cite{flood} and \cite{bhamidi2017universality} by exploring balls of a particular size around the vertices, and constructing the graph at the same time.
The shortest weighted path between two vertices $u$ and $v$ will be described by the first time collision of the two exploration balls around $u$ and $v$. Another way to understand that is to imagine water percolating in the graph started from two different nodes. In this case, the growing exploration ball around a vertex $u$ at a time $t$ can be seen as the set of nodes reached by the flow until this time starting from $u$. \\ We now give a precise definition of this exploration process

\begin{itemize}
\item At time 0, we look at the $d_u$ half-edges incident to $u$ and $d_v$ half-edges incident to $v$ and remove all those forming self-loops at $u$ or $v$. If two half-edges incident to $u$ and $v$ respectively are matched, they form a collision edge and we assign to it a random weight according to $G$. Assign random weights with distribution $G$ for the remaining half-edges and write $A(0)$ for these unmatched half-edges.  
\item Wait until the minimum of lifetimes, denoted by $T_1$, of the active half-edges is reached (the minimum is unique almost surely since $G$ is continuous). 
\item The corresponding half-edge, denoted by $e^*$,  with weight $T_1$ is matched with any other randomly chosen free half-edge, and give weight $T_1$ to the newly formed edge. 
\item Remove the newly discovered half-edges that are part of loops or cycles, update $A(T_1)$ by removing $e^*$ from $A(0)$ and adding the remaining newly discovered free half-edges. 
\end{itemize}
\textit{Remark: This exploration process shows how to explore a neighborhood of a vertex by looking at the random weights on the edges and constructing the graph at the same time by random matching of the half-edges. The order in which we choose the half-edges to be paired in the configuration model does not affect this exploration process. In the sequel of this paper, we will be using different variants of this exploration process that will be useful to get upper and lower bounds for the diameter of the random graph in order to prove Theorem \ref{thm1} given in the next section. }

\section{Main theorem, overview of the approach}
\noindent We now state the main result of this paper:
\thm{\label{thm1}Let $CM_n(\mathbf{d})$ be a random graph constructed according to the configuration model with $i.i.d$ edge weights with common law $G$ satisfying condition \eqref{*} and the degree sequence satisfying \hyperref[enum]{Condition 1}. Then we have

\begin{align*}
\frac{diam_w(CM_n(\mathbf{d}))}{\log n}&\stackrel{p}{\lto} \frac{1}{\alpha}+\frac{2}{cd_{\min}},\ \ \textnormal{and}\\ \frac{flood_w(CM_n(\mathbf{d}))}{\log n}&\stackrel{p}{\lto} \frac{1}{\alpha}+\frac{1}{cd_{\min}},
\end{align*}
where $\alpha$ is the Malthusian parameter of a branching process with degree law $\widehat{\textbf{p}}$ and edge weight distribution $G$ for the particles, see \cite{athreya2004branching}.}

In the penultimate section we establish a \textquote{converse}

\thm{\label{thm2}Let $CM_n(\mathbf{d})$ be a random graph constructed according to the configuration model with $i.i.d$ edge weights with common continuous law $G$ with the degree sequence satisfying \hyperref[enum]{Condition 1}. If we have

\begin{align*}
\frac{diam_w(CM_n(\mathbf{d}))}{\log n}&\stackrel{p}{\lto} \frac{1}{\alpha}+\frac{2}{cd_{\min}}\\ \frac{flood_w(CM_n(\mathbf{d}))}{\log n}&\stackrel{p}{\lto} \frac{1}{\alpha}+\frac{1}{cd_{\min}},
\end{align*}
then \eqref{*} holds with value $c\in (0,\infty)$.}
Theorems \ref{thm1} and \ref{thm2} generalize the result in \cite{flood}. According to these two theorems, the weighted diameter and flooding time on the configuration model are of the order of $\log n$ as $n\to \infty$ if and only if the weight distribution $G$ belongs to a set of light-tailed distributions satisfying \eqref{*}.\\ 
We will focus in this paper on proving Theorems \ref{thm1} and \ref{thm2} for the diameter of the graph. Based on the same techniques, we show in section \ref{flooding} how we can get the desired asymptotics in these theorems for the flooding time.\\
\noindent The idea of the proof, for the diameter, is to study the growth of a ball centered, according to the weighted distance, at a certain vertex, and the time needed until any two such balls intersect. The same tools are used to study the behavior of the flooding, and so its proof is almost straightforward once the result for the diameter is proven. The coupling argument for the growth of the balls and the construction of the graph at the same time are explained in detail in \cite{flood}. The idea is to start from a vertex with a certain number of half-edges $d_i$ and assign to each of them $i.i.d$ weights according to $G$. 
\\ According to \cite{bhamidi2017universality}, we know that the typical size of the balls around two uniformly chosen vertices $u$ and $v$ for collision is of order $\sqrt{n}$. In our case, since we are studying the weighted diameter of the graph and thus considering all the $\binom{n}{2}$ pairs of vertices, we will see that we need to explore the neighborhood of the vertices until a size of the order $\sqrt{n\log n}$.\\ 
\noindent The proof will be divided into two parts. We will first prove, in section \ref{upperbound}, the upper bound for the diameter, by finding first an upper bound for the time needed to reach a size of $K\log n$ half-edges while exploring the neighborhood of a vertex, where $K$ is a constant that is chosen to be large enough and will be useful to prove the upper bound (see Theorem \ref{sn}). We then show that for any $\epsilon > 0$, with high probability, the time needed for all these $K\log n$ half-edges to connect to new vertices is less than $\frac{\log n}{cd_{\min}} (1 + \epsilon) $. 
We then show that we need at most a time $\sqrt{\log n}$ before having at least $K\log n/2$ new splittings, each one of them giving at least 2 new half-edges (so we have at least $K\log n$ new processes). 
Then, using a coupling argument, we show that, as $n\to \infty$, there exist at least two sub-processes, among the $K\log n$ starting subprocesses, that will reach together a size of the order of $\sqrt{n\log n}$ in a time bounded by $(1 + \epsilon ) \frac{1}{2\alpha}\log n$.\\ 
In section \ref{slb}, we show the lower bound for the diameter by finding at least two vertices $u$ and $v$ such that, for any $\epsilon>0$,
\begin{equation*}
dist_w(u,v)\geq \frac{(1-\epsilon)\log n}{\alpha}+\frac{2(1-\epsilon)\log n}{cd_{\min}},\ w.h.p.
\end{equation*}
Finally in the last section, we describe the behavior of the flooding time, as $n\to \infty$, using the same arguments and results as for the diameter.

\section{Upper bound}\label{upperbound}
The purpose of this section is to provide the upper bound for the diameter needed for Theorem \ref{thm1}.  
As with \cite{bhamidi2017universality}, we will see that for two \textquote{typical}  vertices $v$ and $u$, the weighted distance will correspond to two times the time needed for the discovery process for $v$ and $u$ to reach approximately $\sqrt{n\log n} $ half edges.\\  We will write this time as $U_1(u) + U_2(v) + U_3$, where $U_1(u)$ and $U_2(v) $ are the times for the discovery processes for respectively $u$ and $v$ to gain $K \log n$ half-edges and $U_3$ is twice the subsequent \textquote{time} for the two clusters to meet.  
Typically (for any $K$) the values of $U_1$ and $U_2$ are of order $o(\log n)$ and it is $U_3$ (of order $\log n)$  which dominates.  However, we will see that for exceptional \textquote{slow} points $u$ and $v, \ U_1 $ and $U_2$ can be of order $\log n$.  We will also see that for $K$ fixed but large, the term $U_3 /log(n)$ is very close to $1/ \alpha $ uniformly over $u$ and $v$.\\  
In subsections \ref{ctbp} and \ref{klog}, our chief aim is to bound the tails of the random variable $U_1(u)$ \big(or $U_2(v)$\big)  uniformly over all vertices. 
We define for a vertex $v \in V$ and positive $C$, the random variable $T_C (v) \ = \ \inf \big\{t\ \big|\ $ the discovery process for $v$ has at least $C $ half-edges$  \big\}$.   When a vertex $v$ is given or fixed we drop the dependence on $v$ and write $T_C$.  The principle result for this section (which will be proven in the second subsection) is 
\prpo{
\label{proptail}
For any $\epsilon > 0 $ and any $K < \infty $, we have
$$
P \left( \max_{v \in V^n} T_{K\log n} (v) < \frac{(1+\epsilon) \log n}{c d_{\min}}\right) \  \rightarrow \ 1 \mbox{ as } n \ \rightarrow \ \infty .
$$
}
	{\it Remark: Here and elsewhere we write $T_{K\log n}$ and not $T_{\lceil K\log n\rceil}$  where an a priori non integer value is offered for an integer argument.}\\
This result evidently follows immediately from the lemma below which is shown in the second subsection.
\lm{\label{lub}For any $\epsilon>0$, there exists $h>0, \delta > 0$ such that, for sufficiently large $n$, we have 
\begin{equation*}
\pr\left(T_{K\log n} (v)\geq \frac{(1+\epsilon)\log n}{cd_{\min}}\right)\leq n^{-(1+\delta )}h.
\end{equation*}
}
\subsection{The result for tree-branching process}\label{ctbp}
We consider, in this subsection, a continuous time generalized (non Markov) branching process $ (Z(t): t \geq 0)$ with 
$Z(0) = d_{\min} $ and so that individuals have a lifetime distributed independently as $G$ at the end of which they split into a random number of \textquote{offspring} which has law size 
equal to the biased distribution $\{\widehat{p}_k\}_{k \geq d_{\min}-1}$ given in \eqref{bias}.  So by abuse of notation in this subsection $T_{K\log n} $ will denote the time for the branching process to attain population size $K\log n$.  We prove
\lm{
\label{lembranch}
For any $\epsilon > 0$, there exists $C$ and $\delta > 0 $ so that
$$
\pr \left( T_{K\log n}  > \frac{(1+2\epsilon) \log n}{c d_{\min}}\right) \ <  \ \frac{C}{n^{1+\delta}} \mbox{ as } n \ \rightarrow \ \infty .
$$
}
In the following subsection (\ref{klog}), we adapt the approach presented here to show the same result in the general case (Lemma \ref{lub}), where the exploring ball around a vertex up to size $K\log n$ contains cycles, which is the case of any realization of the configuration model w.h.p..

We fix $v$.  To analyze $T_{K\log n}$, which represents the time needed for the continuous time branching process starting from $v$ to reach $K\log n$ half-edges, we use some comparisons with simpler objects.  This is chiefly to deal with the absence of the the memoryless property for general distribution $G$ satisfying \eqref{*}.
For the branching process extra edges can only serve to reduce the random variable $T_{K\log n}$, so we may (and shall) take the number of offspring to be deterministically equal to $d_{\min} -1 \ \geq \ 2$, since we are looking for an upper bound for $T_{K\log n}$.   This being the case we may regard our branching process as derived from a rooted tree where the root has $d_{\min} $ \textquote{offspring}  and subsequent vertices have $d_{\min} -1 $ offspring.  
We associate to each edge $e$ of the tree the random variable $X_e$ where the $X_e$'s are i.i.d. random variables distributed as $G$. The idea is to use condition \eqref{*} in order to stochastically upper bound it by an exponential random variable and use these exponential random variables to find the desired upper bound. 
Our first real comparison process comes by \textquote{freezing} the births of the branching process $Z$ beyond the $(\log n) ^ \gamma $ generation for some fixed $0 < \gamma <1 $. 
 Alternatively we can see this as changing all the variables $X_e $ corresponding to edges from a $(\log n)^ \gamma $ generation vertex to equal infinity. 
  Such a process must necessarily reach (at a random time) the configuration of $d_{\min} (d_{\min} - 1) ^{(\log n)^ \gamma - 1} $ individuals which is bigger than $K\log n$ for large $n$. Thus, writing $T^{'} _{K\log n} $ for the time this modified branching process has $K\log n $ individuals, we obviously have $ T_{K\log n} \ \leq \ T^{'}_{K\log n} $ and so an upper bound on tail probabilities for the latter will serve for the former.  
The next comparison process involves changing the $X_e $ random variables to shifted exponentials:  Property \eqref{*} entails that for each $\epsilon > 0 $ there exists $R_ \epsilon  < \ \infty $ so that 
$$
\forall x \geq R_ \epsilon \ \quad 1-G(x) \leq e^{-c(1- \epsilon ) x},
$$
from which it follows that $G$ is stochastically dominated by the exponential distribution with parameter $c(1- \epsilon ) $ shifted by $R_ \epsilon $ to the right.  
By abuse of notation we write
\begin{equation}\label{gsttt}
  G  \stackrel {st} {\leq } R_ \epsilon +  \mathcal Exp (c(1- \epsilon )).
\end{equation}

Accordingly we can couple random variables $X_e $ with i.i.d.  $ \mathcal Exp (c(1- \epsilon )) $ random variables $X''_e $ so that for each edge $e, \ X_e \leq X''_e  + R_ \epsilon $.

Our final comparison involves $T^{''}_{K\log n} $, the time for the branching process with variables $G^{''}(e) $ to have $K\log n $ individuals where again no birth after generation $\log ^ \gamma (n) $
are permitted.  
$T^{''}_{K\log n} $ is obviously easier to deal with than its preceding objects.  We also note that while in general $T^{''}_{K\log n} $ maybe less than $T^{'}_{K\log n} $, given that we only allow generations up to $\log ^ \gamma (n) $, we have
$$
T^{'}_{K\log n}  \leq  T^{''}_{K\log n} + \log ^ \gamma (n) R_ \epsilon
$$
and that the latter term is negligible compared to $\log n $ as $n$ becomes large. 

So our proof of Lemma \ref{lembranch} has been reduced to proving
\lm{ \label{branchbound}
For $\epsilon > 0, $ there exists $h, \delta > 0 $ so that for all $n$ large
$$
\pr\left(  T^{''}_{K\log n} > \frac{(1+ 2\epsilon)\log n}{ c d_{\min}}\right) \ <\ \frac{h}{n^{1+ \delta}} .
$$
}

Before proving this lemma we will need an elementary counting result for regular trees.  In our deterministic branching model each birth increases $Z$, the population size by $d_{\min}-2$.  Thus it will increase the jump rate of $Z$ by $d_{\min}-2$ unless the $d_{\min}-1$ offspring are of generation $\log ^ \gamma (n) $ in which case the rate is reduced by $1$. \\ Writing $L$ for the number of splittings needed to reach size $K\log n$, we have
\begin{equation}\label{ll}
d_{\min}+L(d_{\min}-2)=K\log n,\ \ \Longrightarrow L=\frac{K\log n-d_{\min}}{d_{\min}-2}.
\end{equation}
As mentioned before, all integer variables used in this paper (that represent a certain number of splittings, generations or number of half-edges....) are written without $\lceil\ \rceil$ brackets for simplicity.\\
We want to find an upper bound for $T^{''}_{K\log n}$ which will also serve as an upper bound (stochastically) for $T_{K\log n}$. To do so, we want to show that, in spite of the restrictions imposed on our modified process (only $d_{\min}$-degree vertices, freezing half-edges at generation $(\log n)^\gamma$), the number of half-edges discovered after each splitting is still sufficiently large in order to reach size $K\log n$ in a time of order $\log n$ at maximum. \\ By \eqref{ll}, the time $T^{''}_{K\log n}$  is equal to the sum of $L=\frac{K\log n -d_{\min}}{d_{\min}-2} $ times between jumps of process $Z$.   That is $T^{''}_{K\log n} \ = \ \sum_{i=1}^{L} F_i$ where $F_i $ is the time between the $i-1 $'th jump and the $ i$'th.  Conditional upon the generational information of the jumps up to the $i-1$'th jump, the random variable $F_i $ is an exponential random variable of parameter $c(1- \epsilon) $ times an integer which is measurable with respect to the information up to the $(i-1)$'th jump.  Up until $ i = \log^ \gamma (n) -1$ the parameter of $F_i $ is nonrandom and equal to $d_{\min} + (i-1)(d_{\min}-2)  $ times $c(1- \epsilon ) $.  
Thereafter the rate can rise or fall.  The lemma below (which is far from optimal but equal to our needs)  records that after this point the parameter of $F_i $, that we denote by $f_i$, has a large lower bound.

\lm{ \label{graphstuff}
For $n$ large and for all $ \log^ \gamma (n) -1 \leq    i \leq L$, $F_i $ has parameter $f_i$ satisfying
\begin{equation*}
f_i\geq [\log^ \gamma (n) / 2  ]\times c(1- \epsilon ).
\end{equation*}
}
\pf{
Let $M$ be the number of splittings at which generation $(\log n)^\gamma$ is  reached for the first time.
Obviously $M \geq  \log^ \gamma (n) -1 $.  If $(\log n)^\gamma-1\leq i<M$, we have
\begin{equation*}
f_i=(d_{\min}+(i-1)(d_{\min}-2))\times c(1-\epsilon)\geq (d_{\min}+((\log n)^\gamma-2)(d_{\min}-2))\times c(1-\epsilon)\geq [\log^ \gamma (n) / 2 ]\times c(1- \epsilon ),
\end{equation*}
for large $n$ and using that $d_{\min}\geq 3$. \\
Suppose now that $i\geq M$. After the $M'th $ jump there is a path  from the root
$v$ to one of the generation $(\log n)^\gamma$  individuals.  Let $u_j$ be the vertex belonging to this path at generation $j\leq (\log n)^\gamma$.
Notice that at least $\log^\gamma (n)/2$ number of generations can be discovered in the subtrees having roots $u_1,u_2,\cdots u_{\frac{ (\log n)^\gamma}{2}}$ before each one of them reaches a total of $\log^\gamma (n)$ number of generations. This means that if one of these subtrees has only free half-edges at generation $\log^\gamma (n)$ (that won't contribute to the jump rate since they are \textquote{freezed}), then their number is at least $(d_{\min}-1)^{\log^\gamma(n)/2}$.\\
But, for $n$ sufficiently large, we have
\begin{equation*}
(d_{\min}-1)^{\frac{(\log n)^\gamma}{2}}>K\log n.
\end{equation*}
This contradicts $i\leq L$ where $L$ is the number of splittings to reach size $K\log n$. This shows that each of these $\frac{\log^\gamma(n)}{2}$ subtrees has at least one free half-edge belonging to one of the first $\log^\gamma(n)-1$ generations of the main branching process. 
In other words, we see that (provided $n$ is large enough) before time $T^{''}_{K\log(n)}$ each one of these subtrees must \textquote{supply} a jump rate of at least $c(1- \epsilon) $ and so we have $f_i\geq \log^\gamma (n)/2\times c(1-\epsilon)$.
}


\begin{myproof}
Using the same notations as in Lemma \ref{graphstuff}, we write\\ $T^{''}_{K\log n}=\mathcal Exp(\lambda d_{\min})+\sum_{i=2}^L F_i:=\mathcal Exp(\lambda d_{\min})+T$ where $\lambda:=c(1-\epsilon)$. We want to show that, for $a:=\frac{(1+ 2\epsilon)\log n}{ cd_{\min}}$ and for any $0<s<a$, that there exist $h,\delta>0$ such that
\begin{equation*}
\pr(T\geq a-s)<he^{\lambda d_{\min}s}n^{-(1+\delta)},
\end{equation*}
for $n$ large. This will finish the proof since
\begin{equation}\label{ptot}
 \pr(T^{''}_{K\log n}\geq a)=e^{-\lambda d_{\min}a}+\int_0^a \lambda d_{\min}\pr(T>a-s)e^{-\lambda d_{\min}s}ds\leq ahn^{-(1+\delta)}+e^{-\lambda d_{\min} a}\sim h'n^{-(1+\delta)},
 \end {equation}
 for a certain $h'>0$.
 Using the Markov inequality and Lemma \ref{graphstuff}, we obtain for $T$,  recalling that $L=\frac{Klog(n) -d_{\min}}{d_{\min}-2} $, 
\begin{align*}
\pr(T\geq a-s)&=\pr\big(e^{\lambda d_{\min}T}\geq e^{\lambda d_{\min} (a-s)}\big)\leq \E\big[e^{\lambda d_{\min} T}\big]e^{-\lambda d_{\min} a}e^{\lambda d_{\min}s}\\
&\leq \prod_{i=2}^{(\log n)^\gamma-2}\left(1+\frac{\lambda d_{\min}}{ ((i-1)(d_{\min}-2)+d_{\min})\lambda -\lambda d_{\min}}\right)\\& \times  \prod_{i=(\log n)^\gamma-1}^{L}\left(1+\frac{\lambda d_{\min}}{(\log^ \gamma (n) / 2 )\times \lambda -\lambda d_{\min}}\right)e^{-\lambda d_{\min}a}e^{\lambda d_{\min}s}\\
&\leq \exp\left[\sum_{i=2}^{(\log n)^\gamma-2}\frac{d_{\min}}{ (i-1)(d_{\min}-2)}\right]\\&\times \exp\left[\sum_{i=(\log n)^\gamma-1}^{L}\frac{d_{\min}}{ (\log^ \gamma (n) / 2 )- d_{\min} }\right]e^{-\lambda d_{\min}a}e^{\lambda d_{\min}s}\\
&\lesssim e^{\frac{d_{\min}}{d_{\min-2}}\log^\gamma (n)}\times e^{\frac{2d_{\min}L}{(\log n)^\gamma(1-o(1))}} \times e^{-\lambda d_{\min}a}e^{\lambda d_{\min}s}\\
&=e^{-(1-\epsilon)(1+2\epsilon)\log n(1-o(1))}e^{\lambda d_{\min}s}
\end{align*}
Therefore, for large $n$ and $\epsilon$ sufficiently small, there exists $\delta>0$ such that
\begin{equation*}
\pr\left(T\geq \frac{(1+ 2\epsilon)\log n}{ cd_{\min}}\right)\leq n^{-(1+\delta)}e^{\lambda d_{\min}s}.
\end{equation*}
This concludes the proof by \eqref{ptot}
 \end{myproof}
\subsection{The result for the general case}\label{klog}
\noindent We showed in the previous section the upper bound for the time needed to reach $K\log n$ half-edges starting from a random vertex, assuming that no cycles or loops occur before that time.  
We show in this section that the same bound holds in case we have one ore more cycles. We say that two paths starting at a vertex $v$ generate a cycle whenever they have another vertex $v'$ in common. We extend this definition to the case where two half-edges incident to the same vertex are matched together and hence forming a loop at this vertex. \\ 
We will first show that, with high probability, we need at maximum $\frac{(1+\epsilon)\log n}{cd_{\min}}$ amount of time, starting from a vertex $v$, to reach $K\log n$ half-edges if we only have exactly 1 cycle in the exploration process. Then, we will show that the probability of having two or more cycles during this process is very small compared to $n^{-1-\delta}$ for $0<\delta<1$, as $n\to \infty$. Hence, this will be sufficient to prove the upper bound of $T_{K\log n}$ in the general case.\\

\subsubsection*{Exactly one cycle:}
Suppose at first that we have exactly one cycle before reaching $K\log n$ half-edges. In this case, the maximal degree of a newly discovered vertex should be less than $K\log n$ (because we stop when we reach $K\log n$ half-edges). On the other hand, we have at maximum $K\log n$ half-edges that can create a cycle during this exploration process (before reaching $K\log n$ half-edges). Hence, the probability of having a cycle can be bounded as follows
\begin{equation}\label{1cycle}
\pr(\textnormal{one cycle at the $i$th splitting})\leq \frac{K\log n\times K\log n}{l_n-i}\lesssim \frac{C(\log n)^2}{n},
\end{equation}
where $C=K^2/m$ and where we used that $\frac{l_n}{n}\to m$ when $n\to \infty$.\\
We will consider the $d_{\min}$-regular case where all newly added vertices have degree $d_{\min}$. Then we will show that even in this case, the time needed to reach size $K\log n$ (even if there are cycles and loops) is upper bounded by $\frac{\log n}{cd_{\min}}$ with high probability. \\ In order to justify the restriction to the $d_{\min}$-regular case, we use a similar comparison argument as in the previous section in order to simplify the current setting.\\ If we have one cycle before reaching $K\log n$ half-edges, we remove the two half-edges that formed a cycle, to obtain an almost $d_{\min}$-regular tree. \\ Let $T'''_{K\log n}$ be the time needed for this almost $d_{\min}$ regular tree to reach size $K\log n$ (by always connecting to new vertices with degree $d_{\min}$). This amount of time is clearly greater or equal to the one in the previous case where only one cycle occurs and no restrictions on the degrees of the vertices are made. \\ Thus, its sufficient to show that Proposition \ref{proptail} also holds for $T'''_{K\log n}$.\\
Notice first, as in the previous section, that the number $S_i$ of alive particles after the $i$th splitting in the $d_{\min}$-regular branching process is given by:
\begin{equation*}
S_i= d_{\min}+(d_{\min}-2)(i-1).
\end{equation*}
After removing the two-half edges that formed a cycle at the $i$th splitting (for a certain integer $i$), there are $d_{\min}+(i-1)(d_{\min}-2)-2$ remaining half-edges. Therefore, we need at most two new splittings to obtain at least $S_i$ half edges since
\begin{equation*}
d_{\min}+(i-1)(d_{\min}-2)-2+2d_{\min}-4\geq S_i,\ \ d_{\min}\geq 3.
\end{equation*}
We let $\tau_1$ be the time spent until the $i$th splitting, $\tau_2$ the time to reach at least $S_i$ again after removing the two bad half-edges and $\tau_3$ be the remaining time to reach $K\log n$ half-edges. We write $R_j=(i-1+j)(d_{\min}-2)$ for $j=1,2$. We obtain, by a simple computation, for $\epsilon>0$,
\begin{align*}
\pr(\tau_2\geq \frac{\epsilon\log n}{cd_{\min}})&=\pr(\mathcal Exp(R_1c(1-\epsilon))+\mathcal Exp(R_2c(1-\epsilon))\geq \frac{\epsilon\log n}{cd_{\min}})\\
&=\frac{R_2c(1-\epsilon)e^{-R_1c(1-\epsilon)\frac{\epsilon\log n}{cd_{\min}}}-R_1c(1-\epsilon)e^{-R_2c(1-\epsilon)\frac{\epsilon\log n}{cd_{\min}}}}{R_2c(1-\epsilon)-R_1c(1-\epsilon)}\\
&\leq \frac{R_2e^{-R_1c(1-\epsilon)\frac{\epsilon\log n}{cd_{\min}}}}{R_2-R_1}
=R_2e^{-R_1c(1-\epsilon)\frac{\epsilon\log n}{cd_{\min}}}\\
&\leq (1+K\log n)(d_{\min}-2)n^{-\frac{d_{\min}-2}{d_{\min}}\epsilon(1-\epsilon)}.
\end{align*}
We write $C_i$ for the event \textquote{Exactly one cycle occurred, at the $i$th splitting}. We finally obtain, using Lemma \ref{branchbound} and \eqref{1cycle},
\begin{align*}
\pr\left(\tau_1+\tau_2+\tau_3\geq\frac{(1+3\epsilon)\log n}{cd_{\min}}, C_i\right)&\leq
\frac{C(\log n)^2}{n}\times\bigg[\pr\left(\tau_1+\tau_2+\tau_3\geq\frac{(1+3\epsilon)\log n}{cd_{\min}},\tau_2\geq\frac{\epsilon\log n}{cd_{\min}}\right)\\&+\pr\left(\tau_1+\tau_2+\tau_3\geq\frac{(1+3\epsilon)\log n}{cd_{\min}},\tau_2\leq\frac{\epsilon\log n}{cd_{\min}}\right)\bigg]\\
&\leq \left(\pr\left(\tau_2\geq\frac{\epsilon\log n}{cd_{\min}}\right)+\pr\left(\tau_1+\tau_3\geq\frac{(1+2\epsilon)\log n}{cd_{\min}}\right)\right)\\&\times\frac{C(\log n)^2}{n}\\
&\leq \left((1+K\log n)(d_{\min}-2)n^{-\frac{d_{\min}-2}{d_{\min}}\epsilon(1-\epsilon)}+hn^{(-1+\epsilon)(1+2\epsilon)}\right)\\&\times\frac{C(\log n)^2}{n}.
\end{align*}
Writing $C'$ for the event \textquote{Exactly one cycle occurred before time $K\log n$}, we get 
\begin{align*}
\pr\left(\tau_1+\tau_2+\tau_3\geq\frac{(1+3\epsilon)\log n}{cd_{\min}}, C'\right)&\leq\left((1+K\log n)(d_{\min}-2)n^{-\frac{d_{\min}-2}{d_{\min}}\epsilon(1-\epsilon)}+hn^{(-1+\epsilon)(1+2\epsilon)}\right)\\&\times\frac{KC(\log n)^3}{n}.
\end{align*}
Hence, taking the union of this event over all the vertices of the graph and writing $h_1$ for this probability, we get
\begin{equation*}
h_1\leq \left((1+K\log n)(d_{\min}-2)n^{-\frac{d_{\min}-2}{d_{\min}}\epsilon(1-\epsilon)}+hn^{(-1+\epsilon)(1+2\epsilon)}\right)\times KC(\log n)^3\to 0, \ n\to \infty.
\end{equation*}
This shows that, starting from any vertex, we need with high probability at most $ \frac{\log n}{cd_{\min}}$ amount of time to reach size $K\log n$ in the exploration process around this vertex, assuming that we have at most one cycle in this exploration process.
\subsubsection*{Two or more cycles:}
On the other hand, using \eqref{1cycle}, the probability $h_2$ of having two or more cycles before reaching size $K\log n$ (starting from a fixed vertex) is bounded, for large $n$ by 
\begin{equation*}
 h_2\lesssim \left(\frac{C(\log n)^2}{n}\right)^2=\frac{C^2(\log n)^4}{n^2}.
\end{equation*}
We thus obtain, by writing $C''$ for the event \textquote{Two or more cycles occurred before time $K\log n$},
\begin{equation*}
\pr\left(\left(\tau_1+\tau_2+\tau_3\geq\frac{(1+\epsilon)\log n}{cd_{\min}}, C''\right)\right)\leq \pr(C'')\leq \frac{C^2(\log n)^4}{n^2}.
\end{equation*}
Hence, taking the union of this event over all the vertices of the graph and writing $h_3$ for this probability, we get
\begin{equation}\label{tmc}
h_3\leq  \frac{C^2(\log n)^4}{n}\to 0,\ n\to \infty.
\end{equation}
\subsection{Time for at least $\frac{K\log n}{2}$ splittings}\label{P2}
Once we reach a number of $K\log n$ half-edges in the exploration process (with a time upper bounded with high probability by $\frac{\log n}{cd_{\min}}$ as seen in the previous section), we denote by $(R_i^v)_{i\leq K\log n}$ the random variables corresponding to the remaining times on the $K\log n$ half-edges obtained in the previous branching process with the root $v$ and we write $(X_i^v)_{i\leq K\log n}$ the corresponding random variables with cdf $G$ representing the total weights on these half-edges. So we have $R_i^v\leq X_i^v$ for all $i\leq K\log n$ and any vertex $v$.\\
In the case of an exponential distribution for the edge-weights (with rate 1 for example), the $R_i^v$'s have also the same exponential distribution by the memorylessness property of the exponential law. In this case, we can study the time until collision between two balls around vertices $u$ and $v$, where both of these vertices have degree $K\log n$. The law of the waiting time before the first splitting in one of these balls is $\mathcal Exp(K\log n)$ by the memorylessness property of $\mathcal Exp(1)$.\\ Since $X_i^v\sim G$ and $G$ doesn't have the memorylessness property, the random variables $R_i^v$ are not distributed according to $G$.\\
To circumvent this problem, we will show that at least $\frac{K\log n}{2}$ of the $K\log n$ half-edges will be connected to new vertices in an amount of time of the order of $\sqrt{\log n}$. Since $d_{\min}\geq 3$, we will get at least $K\log n$ new half-edges and we have again that the lifetime of these new half-edges is distributed according to $G$. Since $\sqrt{\log n}$ is negligible compared to $\log n$, this waiting time to get these $K\log n$ new half-edges will not affect the upper-bound of the diameter which will be shown to be of the order of $\log n$. 
We can put aside the half-edges that don't connect to new vertices in this $\sim \sqrt{\log n}$ amount of time (there's maximum $\frac{K\log n}{2}$ of these half-edges). Hence, by not considering such half-edges, we will need even more time to reach the typical size for collision starting from at least $K\log n$ (newly discovered) half-edges. It's then sufficient to show that the upper bound still hold in this case.\\
We will show the following theorem
\thm{\label{kl2}Consider the $K\log n$ half-edges that were reached by the branching process around $v$ (as in section \ref{ctbp}), with $(R_i^v)_{i\leq K\log n}$ remaining time on these half-edges before they connect to new vertices. Then we have
\begin{equation*}
n\times \pr\left(\textnormal{At least}\ \frac{K\log n}{2}\ \textnormal{of the}\ R_i^v\textnormal{'s}\geq \sqrt{\log n}\right)\to 0,\ n \to \infty.
\end{equation*}
}
\noindent This will show that, starting from any vertex $v$, with high probability, at least $\frac{K\log n}{2}$ alive particles in the corresponding exploration process, will die in the next $\sqrt{\log n}$ units of time giving birth to at least 2 new particles (since $d_{\min}\geq 3$). 
\pf{\noindent To prove this, we notice first that the number of the explored weighted half-edges needed to obtain $K\log n$ alive particles is less than $3K\log n$. To see that, we consider again the worst case where every vertex has degree $d_{\min}$. In this case, we need $ \frac{K\log n-d_{\min}}{d_{\min}-2}$ splittings to reach size $K\log n$. Therefore, the number of weighted half-edges used in this process is given by, for $n$ sufficiently large,
\begin{equation}\label{bound}
d_{\min}+(d_{\min}-2+1) \frac{K\log n-d_{\min}}{d_{\min}-2} \leq K\log n+d_{\min}-1+ \frac{K\log n-d_{\min}}{d_{\min}-2}\leq 3K\log n.
\end{equation}
We let $X_1,\cdots X_{3K\log n}$ be the maximal set of random variables with cdf $G$ that were discovered during the exploration process until reaching size $K\log n$. We let $A$ be the event that at least $\frac{K\log n}{2}$ of the $X_i^v$'s are bigger than $\sqrt{\log n}$. Since $X_i^v\geq R_i^v$ for all $i$ and vertices $v$, it's sufficient to prove that $\pr(A)\to 0$ faster than $\frac{1}{n}$:
\begin{align*}
\pr(A)&=\sum_{r= K\log n/2}^{3K\log n} \binom{3K\log n }{r}\left(\overline{G}\left(\sqrt{\log n}\right)\right)^r\left(1-\overline{G}\left(\sqrt{\log n}\right)\right)^{3K\log n-r}\\
&\lesssim \binom{ 3K\log n}{ 3K\log n/2}\times \left( \frac{5}{2}K\log n +1\right)e^{-c\sqrt{\log n}\times \frac{K\log n}{2}}\\
&\sim 2^{ 3K\log n}\sqrt{\log n}\ e^{-c\sqrt{\log n}\times \frac{K\log n}{2}}\\
&=e^{-c\sqrt{\log n}\times \frac{K\log n}{2}+ 3K\log n\log 2+\frac{1}{2}\log\log n}\to 0,\ n\to \infty,
\end{align*}
where we used Stirling's approximation 
\begin{equation*}
a!\sim \sqrt{2\pi a}\left(\frac{a}{e}\right)^a,\ a\to \infty,
\end{equation*}
to approximate $\binom{3K\log n}{3K\log n/2}$.}
\corl{The result of Theorem \ref{kl2} also holds in the general case, where one or multiple cycles can be created by two or more of the $K\log n$ half-edges that were obtained by exploring the neighborhood of a vertex $v$.}
\pf{The probability of having two or more cycles is negligible as $n\to \infty$ (see \eqref{tmc}).\\
In the case of one cycle, we can show that Theorem \ref{kl2} holds in the exact same way if we replace $\frac{K\log n}{2}$ by $\frac{K\log n}{2}+1$. In other words, with high probability, at least $\frac{K\log n}{2}+1$ half-edges will be matched within a time of the order of $\sqrt{\log n}$. If one of them creates a cycle (by connecting to one of the $K\log n$ half-edges in the exploration process around $v$), then the other $\frac{K\log n}{2}$ half-edges will connect to new vertices with degree bigger than $d_{\min}\geq 3$, and so we reach at least $\frac{K\log n}{2}\times (d_{\min-1})\geq K\log n$ new subprocesses within a time of the order of $\sqrt{n}$.

}


\subsection{Time for collision starting from $K\log n$}\label{upb}
Starting with $K\log n$ newly discovered subprocesses in the exploration process of vertices $u$ and $v$ respectively (as explained in \ref{P2}), we write $S(u,v)$ for the time spent exploring the $2\times K\log n$ processes before the first collision between the two balls.\\
By section \ref{P2}, we may have more than $K\log n$ free half-edges in each ball at this stage, but we consider only $K\log n$ of them that have weights distributed according to $G$ (whereas other half-edges can have remaining lifetimes that are not distributed according to $G$).\\ The time needed for the collision between the two balls of size $K\log n$ is greater than the time needed for collision for the original balls (that can contain more than $K\log n$ half-edges as explained before). Therefore, it's sufficient to upper bound the time needed to have collision between the two balls of size $K\log n$ each. 
We want to show that we need at maximum $\frac{1+\delta}{\alpha}\log n$ amount of time (with high probability) before the collision happens, for $\delta>0$ arbitrary small. The matching among these half-edges is explained in section \ref{expro}.

\noindent The size-biased distribution corresponding to a distribution $(p_k)_{k\geq 0}$ is given by
\begin{equation}\label{abd}
\widehat{p_k}:=\frac{(k+1)p_{k+1}}{m},
\end{equation}
where $m=\sum_r p_r$ and we let $\nu:=\sum_{k}k\widehat{p_k}$. We will use a slightly modified distribution in order to couple each of the $K\log n$ processes with a continuous branching process with a maximal finite degree $\Delta$. For this, given $\epsilon>0$, we define an i.i.d sequence $(Y_k)_{k\geq 0}$ with distribution
\begin{equation}\label{dist}
q_k^n:=\pr(Y_i=k):=\begin{cases}\left(\frac{k+1}{l_n}\sum_{i=1}^n \mathbbm 1_{\{d_i=k+1\}}-\epsilon \right)\vee 0,& 0<k<\epsilon^{\frac{-1}{3}}\\ 1-\sum_{r=2}^{\Delta-1} q_r^n,& k=0,\\
0,& k\geq \epsilon^{\frac{-1}{3}}\end{cases}
\end{equation}
where $\epsilon>0$ is small, $\Delta$ is the maximal degree in this case verifying $\Delta<\epsilon^{\frac{-1}{3}}$ and $l_n$ is the total number of half edges corresponding to the total of $n$ vertices in the graph. 
Similarly, we define, for every $k\geq 0$
\begin{equation}\label{sb}
p_k^n:=\frac{k+1}{l_n}\sum_{i=1}^n \mathbbm 1_{\{d_i=k+1\}}.
\end{equation}

\subsubsection{Coupling the forward degrees}

\noindent We present now a coupling between the forward degrees (the degree minus 1 of a discovered vertex) and a sequence of i.i.d random variables $(Y_i)_{i\geq 1}$ with common law $\textbf{q}$ given in \eqref{dist} (we write $\textbf{q}$ instead of $\textbf{q}^{\textbf{n}}$ for simplicity). \\
We start by showing that the two balls collide with high probability whenever their sizes exceed $C\sqrt{n\log n}$ for some constant $C$. For this, we first define, for a vertex $u$ and time $s>0$,
\begin{equation*}
B_u(s):=\{h\ |\ h\ \textnormal{ free half-edge at time}\ t\ \textnormal{discovered by the exploration process around}\ u\}.
\end{equation*}
\prpo{\label{coll}For any pair of vertices $u,v\in V^n$, we have with high probability
\begin{equation*}
dist_w(u,v)\leq T_{A_n}(u)+T_{A_n}(v),
\end{equation*}
where $A_n:=\sqrt{3mn\log n}$ and $T_{A_n}(u)$ is the time needed for the ball around $u$ to reach a total of $A_n$ half-edges.}
\pf{Fix two vertices $u$ and $v$ and suppose that $B_u(T_{A_n}(u))$ and $B_v(T_{A_n}(v))$ are disjoint. A free half-edge belonging to $B_u(T_{A_n}(u))$ will be matched uniformly at random with another half-edge in the graph. Therefore, the probability that it is not matched with a half-edge in $B_v(T_{A_n}(v))$ is at most
\begin{equation*}
1-\frac{\sqrt{3mn\log n}}{l_n}.
\end{equation*}
Hence, the probability that the two balls do not intersect immediately is upper bounded by
\begin{equation*}
 \left(1-\frac{\sqrt{3mn\log n}}{l_n}\right)^{\sqrt{3mn\log n}}\lesssim e^{-\frac{3mn(\log n)}{nm}}<n^{-2-\delta},
 \end{equation*}
 where we used that $l_n/n\to m$ and fixed $0<\delta<1$. Thus, by summing over all the pairs of vertices $(u,v)$ in the graph, this probability will tend to 0. 
}
\noindent Let $\tilde{p}_k$ be the probability of having $1\leq k<\epsilon^{\frac{-1}{3}}$ children after a splitting in one of the 2 balls before the collision happens, and supposing that we have at maximum one cycle
. This probability depends obviously on the number of already matched half-edges, but this number is upper bounded by $4\sqrt{3mn\log n}$ by Proposition \ref{coll} and similar computation as \eqref{bound}, so we have, for large $n$
\begin{equation}\label{pq}
\tilde{p}_k\gtrsim\frac{\sum_{i=1}^n (k+1)\mathbbm 1_{\{d_i=k+1\}}-4\sqrt{3mn\log n}}{l_n-4\sqrt{3mn\log n}}\gtrsim p_k^n-\epsilon=q_k.
\end{equation}
We now focus on the evolution of the ball around $u$, by looking at the $K\log n$ processes related to this ball. For the $i$th splitting, $i\geq 1$, we let $\tilde{q}_{k,i}$ be the probability of obtaining $k$ children, and none of them belongs to a cycle or a loop, $\epsilon^{-\frac{1}{3}}\geq k\geq 2$. Then we have, for large $n$
\begin{align*}
\tilde{q}_{k,i}\geq \frac{\sum_{r\geq 1}(k+1)\mathbbm 1_{\{d_r=k+1\}}-2\sqrt{3mn\log n}}{l_n-2\sqrt{3mn\log n}}\times \left(1-\frac{\sqrt{3mn\log n}}{l_n-4\sqrt{3mn\log n}}\right)^k&\gtrsim p_k^n \left(1-k\frac{\sqrt{3mn\log n}}{l_n-4\sqrt{3mn\log n}}\right)\\ &\geq q_k.
\end{align*}
 We write $(U_i)_{i\geq 1}$ for a sequence of i.i.d uniform random variables in $(0,1)$. The branching process approximation used in this section is constructed in the following way:
 \begin{itemize}
 \item For the $i$th splitting the $K\log n$ processes related to $u$, if we have $k$ children with $k>\epsilon^{-\frac{1}{3}}$, then we freeze these half-edges and will not be taken into account later on. 
 \item If $k<\epsilon^{-\frac{1}{3}}$, we keep these half-edges if they don't belong to a cycle and if $U_i\leq \frac{q_k}{\tilde{q}_{k,i}}$.
 \end{itemize}
This gives us a coupling between each of the $K\log n$ processes and a continuous branching process with offspring distribution $q$. \\
\textit{Remark:
By \eqref{pq} and the fact that $q_k=0$ for $k>\epsilon^{-\frac{1}{3}}$, we see that the time needed before the collision of the two balls is larger when considering the branching process with offspring distribution $\textbf{q}$. This shows that the bound for this amount of time (before the collision) in the branching process case is sufficient to bound the actual amount of time in the general case.} 

\noindent We let now $Z_t^n$ be the number of alive particles at time $t$ for a continuous branching process with the law for the children given by \eqref{dist}, bounded by $\Delta$, and continuous cumulative distribution $G$ for the edge weights.
We write also $Z_t$ for the number of alive particles at time $t$ for a continuous branching process with the size-biased law for the children and continuous cumulative distribution $G$ for the edge weights. By \cite[p.152]{athreya2004branching}, we know that, in the supercritical case, 
\begin{equation}\label{c'}
\E[Z_t]\sim c'e^{\alpha t},\ \ c'=\frac{\nu-1}{\alpha \nu^2\int_0^\infty ye^{-\alpha y}dG(y)},
\end{equation}
where $\nu$ is the average number of children at each splitting and $\alpha$ is the Malthusian parameter corresponding to the process, which is the unique solution of
\begin{equation*}
\nu\int_0^\infty e^{-\alpha y}dG(y)=1.
\end{equation*}

\lm{\label{aaa}Let $\nu_n$ and $\nu_n^*$ be the expectations corresponding to $p_k^n$ and $q_k^n$ respectively and $\alpha_n$ and $\alpha_n^*$ the corresponding Malthusian parameter and let $\Delta$ denote the maximal degree of the graph. Then we have
\begin{equation*}
\alpha_n-\alpha_n^*\to 0, \epsilon\to 0.
\end{equation*}
}
\pf{We see first that
\begin{equation*}
\nu_n-\nu_n^*=\frac{\epsilon(\epsilon^{-\frac{1}{3}}-1)\epsilon^{-\frac{1}{3}}}{2}+\sum_{k=\epsilon^{-\frac{1}{3}}+1}^\infty kp_k^n\leq \epsilon^{\frac{1}{3}}+\sum_{k=\epsilon^{-\frac{1}{3}}+1}^\infty kp_k^n.
\end{equation*}
Since $\sum_{k=1}^\infty kp_k^n$ converges uniformly by $(c)$ in \hyperref[enum]{Condition 1}, we have that $\nu_n-\nu_n^*\to 0$ when $\epsilon\to 0$. Let $\alpha_n$ and $\alpha_n^*$ be the corresponding Malthusian parameters of $\nu_n$ and $\nu_n^*$, which are the unique respective solutions of 
\begin{equation*}
H(\alpha_n):=\int_0^\infty e^{-\alpha_n y}dG(y)=\frac{1}{\nu_n},\ \ \  H(\alpha_n^*)=\int_0^\infty e^{-\alpha_n^* y}dG(y)=\frac{1}{\nu_n^*}.
\end{equation*}
We see easily that $H$ is differentiable. Using that $\nu_n\to \E[D^*-1]<\infty$, the derivative of $H$ is bounded as follows for sufficiently large $n$,
\begin{equation*}
H'(\alpha_n)=\frac{-1}{\alpha_n}\int_0^\infty \alpha_n ye^{-\alpha_n y}dG(y)=-\frac{1}{\nu_n}\leq -\frac{1}{2\E[D^*-1]}<0.
\end{equation*}
We then obtain, for a certain $\alpha_0\in ]\alpha_n^*,\alpha_n[$,
\begin{equation*}
|H(\alpha_n)-H(\alpha_n^*)|=|H'(\alpha_0)||\alpha_n-\alpha_n^*|\geq \frac{1}{2\E[D^*-1]}|\alpha_n-\alpha_n^*|.
\end{equation*}
Since $\nu_n-\nu_n^*\to 0$ when $\epsilon\to 0$, we have that 
\begin{equation*}
|\alpha_n-\alpha_n^*|\leq |H(\alpha_n)-H(\alpha_n^*)|\times(2\E[D^*-1])\to 0,\ \epsilon\to 0.
\end{equation*}}
\thm{\label{sn}For $u,v\in V^n$, we let $A(u,v):=\bigg\{S(u,v)>\frac{1+\gamma}{\alpha}\log n\bigg\}$ for $\gamma>0$. Then, for $n$ large enough, there exists $\delta>0$ such that 
\begin{equation*}
\pr\left(A(u,v)\right)<n^{-2-\delta}, 
\end{equation*}
where $\alpha$ is the Malthusian parameter is defined in \eqref{malt} and where we recall that $S(u,v)$ is the time spent exploring the $2\times K\log n$ processes before collision. }
\noindent \textit{Remark: We need to mention that condition \eqref{*} on the tail of the distribution $G$ was used in section \ref{upperbound} to upper bound the time for the exploration process to reach size $K\log n$, as well as for the lower bound in section \ref{slb}, but is not used to prove this theorem.}\\
\noindent This will show that $\pr(\cup_{(u,v)\in V^n\times V^n} A(u,v))\to 0,\ n\to \infty$. In other words, with probability that tends to 1, and using the result of the previous section, we need at maximum $\frac{2}{cd_{\min}}\log n + \frac{1+\gamma}{\alpha}\log n$ amount of time before a collision happens between two exploration process around any two uniformly chosen vertices for an arbitrary small $\gamma>0$.

\pf{Denote $Z_t^{*,n}$ the number of alive particles in a continuous branching process with law $G$ for the edges and probability $q_k$ to have $k$ children for every splitting and every $k\geq 1$. Since we have at least $K\log n$ such processes coming from the exploration balls of $u$ and $v$ respectively, we will write, to simplify the notations, these processes as $U_1(t),\cdots, U_{K\log n}(t)$ for those related to $u$ and $V_1(t),\cdots,V_{K\log n}(t)$ for $v$ and $U_i(t),V_j(t)\sim Z_t^{*,n},\ 1\leq i,j\leq K\log n$.
\\ Let $t^*$ be such that $e^{\alpha_n^*t^*}=\sqrt{3mn\log n}$. We first notice that, for any $\epsilon>0$, there exists $n$ sufficiently large such that
\begin{equation*}
t^*=\frac{1}{\alpha_n^*}\log(\sqrt{3mn\log n})=\frac{1}{2\alpha_n^*}\left(\log(3m\log n)+\log n\right)\leq\frac{1}{2\alpha_n^*}\log n(1+\epsilon).
\end{equation*}
We will now show that there exists at least a pair of processes $(U_i(t),V_i(t))$ that collide before time $t^*$.\\
By Proposition \ref{coll}, $(U_i(t),V_i(t))$ will collide with high probability before time $t^*$ whenever 
\begin{equation*}
U_i(t^*), V_i(t^*)>e^{\alpha_n^*t^*}=\sqrt{3mn\log n}.
\end{equation*} 

\noindent Since $Z_te^{-\alpha t}\stackrel{a.s.}{\to} c'W$ and $W$ has a continuous distribution (see \cite{athreya2004branching}), there exists $0<a<1$ such that, for large $t$,
\begin{equation*}
\pr(U_i(t)<e^{\alpha_n^* t})\leq a.
\end{equation*}
From this, we can easily deduce, using again Proposition \ref{coll} that the probability of collision between $U_i(t^*)$ and $V_i(t^*)$ is greater than $(1-a)^2$ for large $n$ and for a certain $0<a<1$. \\
Hence, the probability that none of these pairs of processes $(U_i(t),V_i(t))$ collide before time $t^*$ is upper bounded by
\begin{equation*}
\pr(A(u,v))\leq (1-(1-a)^2)^{K\log n}=e^{K\log n\log(1-(1-a)^2)}=n^{K\log(1-(1-a)^2)}
\end{equation*}
By taking $K$ sufficiently large, we get that this probability is bounded by $n^{-2-\delta}$ for $\delta>0$. 
\noindent By summing over all the pairs of vertices $(u,v)$ in the graph, we can directly conclude that, with high probability, for any pair $(u,v)$, and after reaching size $K\log n$ around these 2 vertices, there will be collision in less than $2t^*=\frac{1}{\alpha_n^*}\log n(1+\epsilon)$ with high probability. 
By Lemma \ref{aaa}, for any $\epsilon>0$, there exists $\gamma>0$ such that
\begin{equation*}
\alpha_n \frac{1+\epsilon}{1+\gamma}\leq \alpha_n^*\leq \alpha_n \frac{1+\epsilon}{1+\gamma/2}.
\end{equation*}
We conclude that we need at most $\frac{1}{\alpha_n}\log n(1+\gamma)$ amount of time, with high probability, to have collision between the two balls once they reach size $K\log n$ each. This finishes the proof since $\gamma$ is arbitrary small and since $\alpha_n\to \alpha$ as $n\to \infty$.
}

\section{Lower bound}\label{slb}
\noindent The goal of this section is to show that, for any $\epsilon>0$, we have with high probability,
\begin{equation*}
\frac{diam(CM_n(\mathbf{d}))}{\log n}\geq \left(\frac{1}{\alpha}+\frac{2}{cd_{\min}}\right)(1-\epsilon),\ n\to \infty.
\end{equation*}
To do this, it's sufficient to show that for any $\epsilon>0$, we can find two vertices $u$ and $v$ in the graph such that 
\begin{equation*}
dist_w(u,v)\geq \frac{(1-\epsilon)\log n}{\alpha}+\frac{2(1-\epsilon)\log n}{cd_{\min}},\ w.h.p.
\end{equation*}
We will only deal with the worst case, where the exploration process starting from any vertex is a branching process.\\

\subsubsection{Coupling the forward degrees}
While exploring the neighborhood of a vertex $u$, we let $\widehat{d}_i$ be the forward degree (the degree minus one) of the discovered vertex at the $i$th splitting. As in \cite{flood}, we set $\beta_n:=3\sqrt{\frac{m}{\nu-1}n\log n}$ and we present a coupling of $(\widehat{d_i})_{i\leq \beta_n}$ with an i.i.d sequence of random variables. We write $\Delta_n$ for the maximum degree in the random graph on $n$ vertices. By writing the order statistics of the degrees as
\begin{equation*}
d_{(1)}\leq \cdots \leq d_{(n)},
\end{equation*}
we write $\overline{m}^{(n)}:=\sum_{i\geq (\beta_n+1)\Delta_n}d_{(i)}^{(n)}$ and we define the size-biased empirical distribution without considering the $(\beta_n+1)\Delta_n-1$ lowest degrees as
 \begin{equation*}
 \overline{\pi}_k^{(n)}:=\frac{\sum_{i\geq (\beta_n+1)\Delta_n} (k+1)\mathbbm 1_{d_{(i)}^{(n)}=k+1}}{\overline{m}^{(n)}}. 
\end{equation*}
By remark \ref{con-d}, we know that $\Delta_n=o(\sqrt{n/\log n})$. We then conclude that $\Delta_n\beta_n=o(n)$.
Hence, it is easy to see that $\overline{\pi}^{(n)}$ tends to the size-biased distribution $\widehat{\textbf{p}}$ defined in \eqref{abd} as $n\to \infty$.\\
The following lemma, proved in \cite{flood}, will be used for the proof of the main result of this section, Proposition \ref{4.3}. 
\lm{\label{ineq}For a randomly chosen vertex $u$ and $i\leq \beta_n$,
\begin{equation*}
\left(\widehat{d}_u (i)\ | \ \widehat{d}_u (1),\cdots,\widehat{d}_u(i-1)\right)\leq_{st} \overline{D}_i^{(n)},
\end{equation*}
where  $\overline{D}_i^{(n)}$ are i.i.d with distribution $\overline{\pi}^{(n)}$.\\ 
}

\noindent For a vertex $u$ and time $t>0$, let $B'(u,t):=\{v\ | \ dist_w(N(u),v)\leq t\}$ where $N(u)$ represents the set of neighbors of $u$ in the graph.  Based on Proposition 4.3 in \cite{flood}, we show  the following proposition
\prpo{\label{4.3}Let $CM_n(\textbf{d})$ denote the random graph constructed with $n$ vertices and a degree sequence $\textbf{d}=(d_i)_{i=1}^n$. Let $t_n=\frac{(1-\epsilon)\log n}{2\alpha}$, where $\alpha$ is the Malthusian parameter corresponding to a branching process with edge weights distribution $G$ and size-biased offspring distribution $\widehat{\textbf{p}}$. For any two uniformly chosen vertices $u,v\in V_{d_{\min}}$, we have, with high probability
\begin{equation*}
B'(u,t_n)\cap B'(v,t_n)=\emptyset.
\end{equation*}
} 
\pf{According to \cite{athreya2004branching}, in the case of a supercritical age-dependent branching process $(Z_t)_{t\geq 0}$, there exists a constant $c'$ such that
\begin{equation}\label{c'} 
\frac{Z_t}{c'e^{\alpha t}}\stackrel{a.s.}{\lto} W,\ \ \E[W]=1.
\end{equation}

\noindent Let $u\in V_{d_{\min}}$.
We consider the worst case for which $B'(u,t)$ is the union of $d_{\min}$ branching processes growing until time $t>0$ and with forward degree $\overline{D}_i^{(n)}$ for the $i$th splitting. We denote these branching processes by $(Z_t^1)_{t\geq 0},\cdots, (Z_t^{d_{\min}})_{t\geq 0}$. Writing $t'_n:=\frac{(1-\epsilon)\log n}{2\alpha_n}$ with $\alpha_n$ the Malthusian parameter corresponding to $\overline{\pi}^{(n)}$ and $G$, we know that $\alpha_n\to \alpha$ as $n\to \infty$. Let $z_n:=\sqrt{\frac{n}{\log n}}$, we define
\begin{equation*}
q_n:=\pr(Z_{t'_n}^1,\cdots,Z_{t'_n}^{d_{\min}}\leq z_n).
\end{equation*}
Using \eqref{c'}, we have, for any $1\leq d_{\min}$,
\begin{align*}
\pr(Z_{t'_n}^i\leq z_n)\sim \pr\bigg(W\leq \frac{z_n}{c'e^{\alpha_n t'_n}}\bigg)=\pr\left(W\leq \frac{n^\epsilon}{c'\sqrt{\log n}}\right)\to 1,\ n\to \infty.
\end{align*}
This implies that $q_n\to 1$ as $n\to \infty$. Therefore, with high probability, the size of $B'(u,t'_n)$ is bounded by $z_n$. Consequently, the probability of getting a collision edge between $B'(u,t'_n)$ and $B'(v,t'_n)$ is bounded by 
\begin{equation*}
\frac{z_n^2}{l_n}\sim \frac{z_n^2}{nm}\to 0,\ n\to \infty,
\end{equation*}
which completes the proof.}
\textit{\noindent Remark: By \cite{flood}, we have that the number of free half-edges after $\beta_n$ splittings, $S_{\beta_n}(u)$, in the exploration process around a vertex $u$ satisfies for large $n$,
\begin{equation*}
S_{\beta_n}(u)\geq \sqrt{\nu-1}\beta_n\geq \sqrt{3mn\log n}\ \ \textnormal{with probability}\ \geq 1-o(n^{-3/2}).
\end{equation*}
This means that, for any uniformly chosen vertex, we need with high probability at maximum $\beta_n$ splittings before reaching size $\sqrt{3mn\log n}$ which is the typical size order for collision according to Proposition \ref{coll}. Hence, coupling the first $\beta_n$ forward degrees in the exploration process of a given vertex before collision (with another ball) is sufficient with high probability.}\\
\newline \noindent Let $V_{d_{\min}}$ be the set of vertices of degree $d_{\min}$ and let $s_n:=\frac{1-\epsilon}{cd_{\min}}\log n$. A vertex in $V_{d_{\min}}$ is called \textit{bad} if the weights on its $d_{\min}$ connected edges are all greater than $s_n$. We also write $A_u$ for the event that $u$ is a bad vertex. \\
The following lemma shows that the average number of bad vertices in the graph tends to infinity as $n\to \infty$ but is negligible compared to $n$:
\lm{\label{eeee2}For any $\epsilon>0$, there exist $a_\epsilon,b_\epsilon>0$ such that
\begin{equation*}
a_\epsilon p_{d_{\min}}(1+o(1))n^{\epsilon^2}\leq \E[Y]\leq b_{\epsilon}p_{d_{\min}}(1+o(1))n^{2\epsilon}.
\end{equation*}
}
\pf{By condition \eqref{*}, for any $\epsilon>0$, there exist $R'_\epsilon$ such that
\begin{equation*}
G\stackrel{st}{\geq}\mathcal Exp(c(1+\epsilon))-R'_\epsilon.
\end{equation*}
Using this, and writing $X_1,\cdots X_{d_{\min}}$ for the random weights on the half-edges connected to a vertex $u\in V_{d_{\min}}$ , we have
\begin{align*}
\pr(A_u)&=\pr\left(X_1\geq s_n,\cdots X_{d_{\min}}\geq s_n\right)\geq \pr\left(\mathcal Exp(c(1+\epsilon))-R'_{\epsilon}\geq s_n\right)^{d_{\min}}\\&=e^{-c(1+\epsilon)R'_{\epsilon}}\times e^{-c(1+\epsilon)s_n d_{\min}}=a_\epsilon n^{-(1-\epsilon^2)},
\end{align*}
where $a_\epsilon:=e^{-c(1+\epsilon)R'_{\epsilon}}$. From this, we get
\begin{equation*}
\E[Y]=\sum_{u\in V_{d_{\min}}}\pr(A_u)\geq a_\epsilon p_{d_{\min}}(1+o(1))n^{\epsilon^2}.
\end{equation*}
The upper bound for $\E[Y]$ follows similarly using \eqref{gsttt}.
}
\lm{\label{yey}Let $Y=\sum_u \mathbbm 1_{A_u}$ the number of bad vertices in the graph. Then we have
\begin{equation*}
Y\geq \frac{2}{3}\E[Y]\ w.h.p. 
\end{equation*}}
\pf{ 

Using that $\Cov(\mathbbm 1_{A_u},\mathbbm 1_{A_v})$ and $\Var(\mathbbm 1_{A_u})$ are both upper bounded by $\pr(A_u)$, we get,
\begin{align*}
\Var(Y)&=\sum_{u\in V_{d_{\min}}}\Var(\mathbbm 1_{A_u})+\sum_{u\in V_{d_{\min}}}\sum_{v\sim u} \Cov(\mathbbm 1_{A_u},\mathbbm 1_{A_v})\\
&\leq \sum_{u\in V_{d_{\min}}}\pr(A_u)+\sum_{u\in V_{d_{\min}}}\sum_{v\in N(u)} \pr(A_u)
= \E[Y]+\sum_{v\in N(u)} \E[Y]\\
&=(d_{\min}+1)\E[Y].
\end{align*}
By Chebychev's inequality, we obtain, for $A>0$
\begin{align*}
\pr(Y\leq \E[Y]-A)&\leq \frac{\Var(Y)}{A^2}\leq \frac{ (d_{\min}+1)\E[Y]}{A^2}.
\end{align*}
Taking $A=\frac{1}{3}\E[Y]$, we get for large $n$
\begin{equation}\label{y23}
Y\geq \frac{2}{3}\E[Y]\ w.h.p.
\end{equation}

}

\noindent We let $Y'$ denote the number of bad vertices belonging to $B'(a,s_n+\frac{(1-\epsilon)\log n}{\alpha})$ for a uniformly chosen vertex $a$. By Proposition \ref{4.3}, we have, for any vertex $i$,
\begin{equation*}
\pr(A_i,B'(a,t_n)\cap B'(i,t_n)\neq \emptyset)=o(\pr(A_i))\impl \E[Y']=o(\E[Y]).
\end{equation*}
We deduce, by Markov's inequality, that $Y'\leq \frac{1}{3}\E[Y]$ with high probability and thus $Y-Y'>0$ with high probability using Lemma \ref{eeee2}.\\
We write $R= \binom{Y}{2}$ for the number of pairs of distinct bad vertices and $R'$ for the number of pairs of distinct bad vertices at distance at most $2s_n+\frac{1-\epsilon}{\alpha}\log n$.
By Proposition \ref{4.3}, it's easy to see that
\begin{equation*}
\pr(A_u,A_v, B'(u,t_n)\cap B'(v,t_n)\neq \emptyset)=o(\pr(A_u,A_v)).
\end{equation*}
Using this, we get
\begin{equation*}
\E[R']=o(\E[Y^2]).
\end{equation*}
Therefore, with high probability, the difference $R-R'$ is strictly positive by Lemma \ref{eeee2}. We deduce that for any $\epsilon>0$, we can find two vertices that are at distance bigger than $2s_n+\frac{1-\epsilon}{\alpha}\log n$. In other words, we obtain
\begin{equation*}
diam(CM_n(\textbf{d}))\geq 2s_n+\frac{1-\epsilon}{\alpha}\log n.
\end{equation*}
Since $\epsilon$ is arbitrary, this proves the lower bound of the diameter, and thus, by section \ref{upperbound}, we finally obtain
\begin{equation*}
\frac{diam(CM_n(\textbf{d}))}{\log n}\stackrel{p}{\longrightarrow} \frac{1}{\alpha}+\frac{2}{cd_{\min}}.
\end{equation*}

\section{Proof of the converse theorem}\label{equiv}

\noindent A first step to proving Theorem \ref{thm2} is the following which amounts to saying simply that exponential tails are required for the diameter (or flood ) to scale as $\log n $ in the sense of the theorem. 
\lm{  \label{lem60}
If 
$$
  \liminf_{x \rightarrow \infty}  \frac{- \log (\overline G(x))}{x} =  0,
$$
then for every $M \ < \ \infty$,
$$
 \limsup_{n \rightarrow \infty}  \pr\left(diam(CM_n(\textbf{d})) > M \log n\right)  =1
$$
}
{\it Remark: The claimed conclusion obviously contradicts the hypotheses of Theorem \ref{thm2} and so in particular any $G$ for which the hypotheses of Theorem \ref{thm2} hold must possess all moments.}

\begin{proof}
By hypothesis ( it is easily seen)  for every $ \epsilon > 0 , \exists $ a sequence of integers $n_j $ tending to infinity so that 
$$
\forall j \ \frac{-\log\overline{G}(\log n_j)}{\log n_j}  <  \epsilon .
$$
Thus we easily have that with probability tending to 1 as $j \rightarrow \ \infty $ there exist vertices $v \ \in \ V_{dmin} \subset \ V^{n_j} $ so that
$$
\min_{u \sim v } G(u,v) > \ \frac{1}{(d_{\min} + 1) \epsilon} \log n_j
$$
This implies that the diameter or flood for the graph $CM_{n_j}(\textbf{d}) $ must exceed $\frac{1}{(d_{\min} + 1) \epsilon} \log n_j$.
The conclusion follows from the arbitrariness of $\epsilon > 0 $.
\end{proof}

\noindent As usual we establish convergence by suitably bounding the $\limsup$ above and the $\liminf$ below: Theorem \ref{thm2} follows from the two lemmas below.

\lm{ \label{lem61}
For distribution $G$ satisfying the hypotheses of Theorem \ref{thm2}
$$
 \liminf_{x \rightarrow \infty}  \frac{- \log (\overline G(x))}{x} \geq   c
$$
}

\lm{ \label{lem62}
For distribution $G$ satisfying the hypotheses of Theorem \ref{thm2}
$$
 \limsup_{x \rightarrow \infty}  \frac{- \log (\overline G(x))}{x}  \leq  c
$$
}

\begin{myproof61}
Suppose not.  Then there exists $\epsilon > 0 $ and a sequence of integers $n_j $ tending to infinity so that
$$
  \forall j \ \frac{-\log\overline{G}(\log n_j)}{\log n_j} \ < \ c(1-\epsilon ).
$$
We can now argue as in section \ref{slb}.   For random graph $CM_{n_j}(\textbf{d})$ we have that $M_j $ the number of vertices $v$ in $V^{n_j}_{d_{\min}} $ so that 
$\min_{u \in N(v)} G(u,v) \ \geq \ \frac{log(n_j)}{d_{\min}c(1- \epsilon / 2)}$ will satisfy with probability tending to one as $j $ tends to infinity the following two conditions
\begin{enumerate}[(i)]
\item $M_j \ \geq \ c n_j^{1 - (1- \epsilon)/(1- \epsilon /2 )}$ for some universal strictly positive $c$ 
\item For each $ \delta > 0 $ with probability tending to one as $j $ tends to infinity for $u$ and $v$ two randomly chosen vertices among the $M_j $ such vertices
\begin{equation*}
B\left(u, \log(n) \frac{1 - \delta }{\alpha}\right) \cap  B\left(u, \log(n) \frac{1 - \delta }{\alpha}\right)  \ = \ \emptyset.
\end{equation*}
\end{enumerate}
Taking $\delta $ sufficiently small with respect to $\epsilon $ gives 
$$
 \limsup_{n \rightarrow \infty}  \pr\left(diam\left(CM_{n_j}(\textbf{d})\right) > \left(\frac{2}{c d_{\min}(1- \epsilon / 3 )} \ + \ \frac{1}{\alpha }\right)\log n\right)  = 1
$$
which contradicts the hypotheses of Theorem \ref{thm2}
\end{myproof61}

\begin{myproof62}
Suppose not.  In this case there exists $\epsilon > 0 $ and a sequence of integers $n_j $ tending to infinity so that
$$
 \forall j \ \frac{-\log\overline{G}(\log n_j)}{\log n_j}  > c(1+\epsilon ).
$$
We may assume by Lemma \ref{lem61} that 
$$
 \liminf_{x \rightarrow \infty}  \frac{- \log (\overline G(x))}{x} \geq   c
$$
and from this, we can apply the argument of Proposition \ref{proptail}
 and see that as $j $ tends to infinity
$$
\pr\left(  \sup_{v \in V^{n_j}}  T_{K \log n_j} (v) \geq c\log n_j / (1+ \epsilon / 2 )   \right) 
$$
tends to zero.

$$
 \limsup_{n \rightarrow \infty}  \pr\left(diam\left(CM_{n_j}(\textbf{d})\right) < \left(\frac{2}{c d_{\min}(1+ \epsilon / 4 )} \ + \ \frac{1}{\alpha }\right) \log n\right)=1
$$
which again contradicts the hypotheses of Theorem \ref{thm2}.
\end{myproof62}

\section{Flooding}\label{flooding}
\noindent We show in this section, based on the proofs and results obtained in sections \ref{upperbound} and \ref{slb}, that with high probability, the weighted flooding time behaves like $(\frac{1}{\alpha}+\frac{1}{d_{\min}})\log n$ as $n\to \infty$.\\
\begin{itemize}\item We show first that $flood(G)\leq (\frac{1}{\alpha}+\frac{1+2\epsilon}{d_{\min}})\log n$ with high probability. as $n\to \infty$. We let $T_{u,\sqrt{n\log n}}(G)$ be the time needed, starting from vertex $u$, to reach $\sqrt{n\log n}$ half-edges, given that the edge weights have a cdf $G$. We have already shown, in section \ref{upperbound}, that with high probability, for any $v$ vertex of the graph, that $T_{v,\sqrt{n\log n}}(G)\leq(\frac{1}{2\alpha}+\frac{1+2\epsilon}{d_{\min}})\log n$. Hence, it's sufficient to show that $T_{a,\sqrt{n\log n}}(G)\leq (\frac{1}{2\alpha})\log n$ for a randomly chosen vertex $a$.\\ Using similar computations as in Lemma \ref{lub}, we have, for any $\epsilon'>0$,
\begin{equation*}
\pr\left(T_{a,K\log n}(G)\geq \frac{\epsilon'}{2\alpha} \log n\right)\lesssim n^{\frac{-\epsilon'c(1-\epsilon')d_{\min}}{2\alpha}}\to 0,\ \ n\to \infty.
\end{equation*}
By section \ref{upb}, we have that with high probability, the time needed to reach $\sqrt{n\log n}$ half-edges starting from $K\log n$ is smaller than $\frac{1}{2\alpha}\log n$. Therefore, we obtain
\begin{equation*}
\pr\left(T_{a,\sqrt{n\log n}}\geq \frac{(1+\epsilon')\log n}{2\alpha}\right)\to 0,\ \ n\to \infty.
\end{equation*}
Since $\epsilon'$ is arbitrary, we finally obtain
\begin{equation*}
flood(G)\leq \left(\frac{1}{\alpha}+\frac{1+2\epsilon}{d_{\min}}\right)\log n,\ w.h.p.
\end{equation*}
\item For the lower bound, we recall the same notations introduced in section \ref{slb}. Since $Y'\leq \frac{1}{3}\E[Y]$ with high probability and using Lemma \ref{yey}, we have
\begin{equation*}
Y-Y'>0\  w.h.p.
\end{equation*}
In other words, with high probability, there exists a vertex $w$ that does not belong to $B'(a,s_n+\frac{(1-\epsilon)\log n}{\alpha})$ where $s_n:=\frac{1-\epsilon}{cd_{\min}}\log n$. This is equivalent to
\begin{equation*}
flood(G)\geq s_n+\frac{(1-\epsilon)\log n}{\alpha}=\frac{1-\epsilon}{cd_{\min}}\log n+\frac{(1-\epsilon)\log n}{\alpha},\ w.h.p.
\end{equation*}
\end{itemize}

\cleardoublepage
\phantomsection
\bibliographystyle{apt}
\bibliography{bibliography}
\addcontentsline{toc}{chapter}{Bibliography}

\end{document}